\documentclass[12pt]{amsart}
\usepackage{amssymb,latexsym}
\usepackage{array}
\usepackage{enumerate}

\begin{document}

\baselineskip=17pt

\title[Equivariant $K$-theory of flag varieties revisited]{Equivariant
  $K$-theory of flag varieties revisited and related results}
\author{V.Uma} 
\address{Department of Mathematics, IIT Madras,
  Chennai, India} 
\email{vuma@iitm.ac.in} 
\date{}

    

\begin{abstract}
In this article we obtain many results on the multiplicative
structure constants of the $T$-equivariant Grothendieck ring
$K_{T}(G/B)$
of the flag variety $G/B$. We do this by lifting the classes of
structure sheaves of Schubert varieties in $K_{T}(G/B)$ to
$R(T)\otimes R(T)$, where $R(T)$ denotes the representation
ring of
the torus $T$. We further apply our results to describe the
multiplicative structure constants of $K(X)_{\mathbb{Q}}$ where
$X$
denotes the wonderful compactification of the adjoint group of
$G$, in
terms of the structure constants of Schubert varieties in the
Grothendieck ring of $G/B$.
\end{abstract}

\maketitle

\thispagestyle{empty} 

\frenchspacing

\textwidth=13.5cm
\textheight=23cm
\parindent=16pt
\oddsidemargin=-0.5cm
\evensidemargin=-0.5cm
\topmargin=-0.5cm

\def\theequation {\arabic{section}.\arabic{equation}}

\newcommand{\codim}{\mbox{{\rm codim}$\,$}}
\newcommand{\stab}{\mbox{{\rm stab}$\,$}}
\newcommand{\lr}{\mbox{$\longrightarrow$}}

\newcommand{\ra}{\rightarrow}
\newcommand{\blr}{\Big \longrightarrow}
\newcommand{\da}{\Big \downarrow}
\newcommand{\ua}{\Big \uparrow}
\newcommand{\hra}{\mbox{{$\hookrightarrow$}}}
\newcommand{\rt}{\mbox{\Large{$\rightarrowtail$}}}
\newcommand{\dua}{\begin{array}[t]{c}
\Big\uparrow \\ [-4mm]
\scriptscriptstyle \wedge \end{array}}

\newcommand{\be}{\begin{equation}}
\newcommand{\ee}{\end{equation}}

\newtheorem{guess}{Theorem}[section]
\newcommand{\bth}{\begin{guess}$\!\!\!${\bf}~}
\newcommand{\eeth}{\end{guess}}
\renewcommand{\bar}{\overline}
\newtheorem{propo}[guess]{Proposition}
\newcommand{\bpropo}{\begin{propo}$\!\!\!${\bf}~}
\newcommand{\epropo}{\end{propo}}

\newtheorem{lema}[guess]{Lemma}
\newcommand{\blem}{\begin{lema}$\!\!\!${\bf }~}
\newcommand{\elem}{\end{lema}}

\newtheorem{defe}[guess]{Definition}
\newcommand{\bdefe}{\begin{defe}$\!\!\!${\bf}~}
\newcommand{\edefe}{\end{defe}}

\newtheorem{coro}[guess]{Corollary}
\newcommand{\bcor}{\begin{coro}$\!\!\!${\bf}~}
\newcommand{\ecor}{\end{coro}}

\newtheorem{rema}[guess]{Remark}
\newcommand{\brem}{\begin{rema}$\!\!\!${\bf}~\rm}
\newcommand{\erem}{\end{rema}}

\newtheorem{exam}[guess]{Example}
\newcommand{\beg}{\begin{exam}$\!\!\!${\bf}~\rm}
\newcommand{\eeg}{\end{exam}}

\newtheorem{notn}[guess]{Notation}
\newcommand{\bnot}{\begin{notn}$\!\!\!${\bf}~\rm}
\newcommand{\enot}{\end{notn}}

\newcommand{\ctext}[1]{\makebox(0,0){#1}}
\setlength{\unitlength}{0.1mm}
\newcommand{\cl}{{\mathcal L}}
\newcommand{\cp}{{\mathcal P}}
\newcommand{\ci}{{\mathcal I}}
\newcommand{\bz}{\mathbb{Z}}
\newcommand{\cv}{{\mathcal V}}
\newcommand{\ce}{{\mathcal E}}
\newcommand{\ck}{{\mathcal K}}
\newcommand{\cz}{{\mathcal Z}}
\newcommand{\cg}{{\mathcal G}}
\newcommand{\cj}{{\mathcal J}}
\newcommand{\cc}{{\mathcal C}}
\newcommand{\ca}{{\mathcal A}}
\newcommand{\cb}{{\mathcal B}}
\newcommand{\cx}{{\mathcal X}}
\newcommand{\cR}{{\mathcal R}}
\newcommand{\cs}{{\mathcal S}}
\newcommand{\ch}{{\mathcal H}}
\newcommand{\cf}{{\mathcal F}}
\newcommand{\cd}{{\mathcal D}}
\newcommand{\bq}{\mathbb{Q}}
\newcommand{\bt}{\mathbb{T}}
\newcommand{\bh}{\mathbb{H}}
\newcommand{\br}{\mathbb{R}}
\newcommand{\bl}{\mathbf{L}}
\newcommand{\wt}{\widetilde}
\newcommand{\im}{{\rm Im}\,}
\newcommand{\bc}{\mathbb{C}}
\newcommand{\bp}{\mathbb{P}}
\newcommand{\ba}{\mathbb{A}}
\newcommand{\spin}{{\rm Spin}\,}
\newcommand{\ds}{\displaystyle}
\newcommand{\tor}{{\rm Tor}\,}
\newcommand{\bff}{{\bf F}}
\newcommand{\bs}{\mathbb{S}}
\def\ns{\mathop{\lr}}
\def\nssup{\mathop{\lr\,sup}}
\def\nsinf{\mathop{\lr\,inf}}
\renewcommand{\phi}{\varphi}
\newcommand{\co}{{\mathcal O}}
\newcommand{\tT}{{\widetilde{T}}}
\newcommand{\tG}{{\widetilde{G}}}
\newcommand{\tB}{{\widetilde{B}}}
\noindent

\section{Introduction}

Let $G$ be a semi-simple simply connected algebraic group over
an
algebraically closed field $k$. Let $B$ be a Borel subgroup and
$T\subset B$ be a maximal torus.

In this article we construct explicit lifts of the classes of
the structure sheaves of Schubert basis in $K_{T}(G/B)$ in the
ring
$R(T)\otimes_{\bz} R(T)$. For this, we apply techniques similar
to
those developed in the paper of Marlin (see \cite{m}) by
exploiting
the properties of  Demazure
 operators (\cite{d}).

Using these lifts we also give new methods to describe the
multiplicative structure of $K_{T}(G/B)$. More precisely, in
\S2and
\S3, we give closed formulas for the multiplicative structure
constants and also recover some known results on these
constants in
this setting.

This was inspired by the results of Hiller (see Chapter IV of
\cite{hi}) who constructs a basis for $Sym(X^*(T))$ as
$Sym(X^*(T))^{W}$-module by lifting the fundamental classes of
the
Schubert varieties in the cohomology ring $H^*(G/B)$. He
further uses this
basis to develop an algebraic approach to Schubert calculus in
$H^*(G/B)$.

Let $X$ denote the wonderful compactification of the semi-simple
adjoint group $G_{ad}=G/Z(G)$.  Recall that in the main result of
\cite{u}, the images in $K(G/B)$ under $c_{K}$ of the Steinberg basis
of $R(T)$ as an $R(G)$-module are used to describe the multiplicative
structure of the Grothendieck ring $K(X)$ of $X$. Moreover, the
multiplicative structure constants of $K(X)$ as a $K(G/B)$-algebra
involved the images of the structure constants of the Steinberg basis,
which do not have known direct geometric or representation theoretic
interpretations (see Theorem 3.8 and Theorem 3.12 of \cite{u}).

In \S4 we show that the above constructed lifts of the Schubert basis
in $K(G/B)$ to $R(T)$ form a basis of $R(T)_{\mathfrak p}$ over
$R(G)_{\mathfrak p}$, where $\mathfrak{p}$ denotes the kernel of the
augmentation map $R(G)\ra \bz$. These are a new set of bases for
$R(T)_{\mathfrak p}$ as $R(G)_{\mathfrak p}$-module different from the
basis obtained by localization from that defined by Steinberg in
\cite{st}. We then reformulate the results in \S3 of \cite{u} using
these bases instead of the Steinberg bases. Using this reformulation,
in the main result, Theorem \ref{main}, of this paper we prove that
$K(X)_{\bq}$ is a free $K(G/B)$-module generated by classes of the
structure sheaves of Schubert varieties in all $K(G/P)_{\bq}$, where
$P\supseteq B$ is a parabolic subgroup. In particular, we express the
multiplicative structure constants of $K(X)_{\bq}$ as $K(G/B)$-algebra
in terms of the structure constants of the classes of structure
sheaves of Schubert varieties in the Grothendieck ring of flag
varieties.

Thus, although we seem to lose to some extent by going to
rational
coefficients, we do obtain better interpretations of the basis
and the
multiplicative structure of the $K(X)_{\bq}$ by relating it to
the
Schubert calculus in the Grothendieck ring of flag varieties.

\subsection{Notations and Conventions}

As in the introduction, let $G$ be a simply connected
semi-simple
algebraic group over an algebraically closed field $k$, let $B$
be a
Borel subgroup of $G$ and let $T\subseteq B$ be a maximal torus.Let $\Lambda=X^*(T)$ denote the weight lattice. Let $\Phi$
denote the
root system and $\Delta$ the set of simple roots relative to
$B$. Let
$W=N(T)/T$ be the Weyl group of the root system $\Phi$. Let
$B^{-}=w_0Bw_0$ be the opposite Borel subgroup to $B$ where
$w_0$ is
the unique maximal element of the Bruhat order on $W$. Let
$\rho:=\frac{1}{2}\sum_{\alpha\in\Phi^{+}}\alpha$. Let
$\omega_{\alpha}$ denote the fundamental weight corresponding
to the
simple root $\alpha\in \Delta$.

For $w\in W$, let $X_w$ denote the Schubert variety which is the
closure of the Schubert cell $BwB/B$ in $G/B$ and let $X^{w}$ denote
the opposite Schubert variety which is the closure of the opposite
Schubert cell $B^{-}wB/B$ in $G/B$. Thus we have: $X^{w}=w_0X_{w_0w}$.

Let $v,w\in W$. Recall that the {\em Bruhat order} is the order
on $W$
defined by the condition $v\preceq w$ if and only if
$X_v\subseteq
X_w$. Further, $X^v \bigcap X_w$ is nonempty if and only if
$v\preceq
w$; then $X^{v}\bigcap X_{w}$ is a variety called the
Richardson variety and is denoted by $X^{v}_{w}$. Moreover,
$X^{v}_w$ has
two
kinds of boundaries, namely $(\partial X_{w})^{v}:=(\partial
X_{w})\bigcap X^{v}$ and $(\partial X^{v})_{w}:=(\partial
X^{v})\bigcap X_{w}$. Here $\partial X_{w}=\bigcup_{w'\prec
w}X_{w'}$
is the boundary of the Schubert variety $X_{w}$ and $\partial
X^{v}=\bigcup_{v\prec v'}X^{v'}$ is the boundary of the
opposite Schubert variety $X^{v}$. Thus we have $(\partial
X_{w})^{v}=\bigcup_{w'\prec w}X^{v}_{w'}$ and $(\partial
X^{v})_{w}=\bigcup_{v\prec v'}X_{w}^{v'}$ (see Prop.1.3.2 and
\S4.2 of
\cite{br4}).

For $X$ any smooth $G$-variety, let $K_{G}(X)$ denote the
Grothendieck
ring of $G$-equivariant coherent sheaves (or equivalently,
vector
bundles) on $X$.  We have the canonical forgetful homomorphism
$K_{G}(X)\ra K_{T}(X)$. In particular, $R(G):=K_{G}(pt)$ is
the Grothendieck ring of $k$-representations of $G$. Since $G$ is simply
connected, we can identify $R(G)=\bz[\Lambda]^{W}$ via
restriction to
$T$. Furthermore, the structure morphism $X\ra k$ induces a
canonical
$R(G)$-module structure on $K_{G}(X)$. Also, $K(X)$ denotes the
Grothendieck ring of coherent sheaves on $X$ and we have
canonical
forgetful homomorphism $K_{G}(X)\ra K(X)$.

For $w\in W$, let $[\co_{X_w}]_{T}$ (resp. $[\co_{X^w}]_{T}$)
denote
the class of the structure sheaf of the Schubert variety
(resp. opposite Schubert variety) in $K_{T}(G/B)$. Further,
note that
we have the identification $[\co_{X^w}]_{T}=w_0\cdot
[\co_{X_{w_0w}}]_{T}$ in $K_{T}(G/B)$. Recall from \cite{kk}
that the
Schubert classes $\{[\co_{X^{w}}]_{T}\}_{w\in W}$ form a basis
of
$K_{T}(G/B)$ as an $R(T)$-module.

For $\lambda\in \Lambda$, let
$\cl(\lambda):=(\bc_{\lambda}\times
G)/B$, where $B$ acts diagonally, and the $B$ action on the one
dimensional vector space $\bc_{\lambda}$ is given by the
surjection
$B\ra T$ followed by $\lambda$. Then $\cl(\lambda)$ is a
$T$-linearized line bundle on $G/B$ associated to $\lambda$.
Let $\cl^{\lambda}$ denote the class of $\cl(\lambda)$ in
$K_{T}(G/B)$.
Further, we shall denote by $e^{\lambda}$ the class of the
trivial
bundle in $K_{T}(G/B)$ with $T$-action given by $\lambda$.

Let $c^{T}_{K}:\bz[\Lambda]=R(T)\ra K_T(G/B)$ denote the
characteristic homomorphism which sends $e^{\lambda}\in R(T)$
to $\cl^{\lambda}\in K_{T}(G/B)$.

Let $*$ denote the canonical involution in $K_T(G/B)$ defined
by duality of a $T$-vector bundle. This is compatible with the
involution
in $R(T)$ defined by $e^{\lambda}\mapsto e^{-\lambda}$. In
particular,
$*c^{T}_{K}(e^{\lambda})=[(\bc_{-\lambda}\times G)/B]\in
K_{T}(G/B)$. 

If $Y\supset Z$ are closed $T$-stable subvarieties of a $T$-variety
$X$, then $\co_{Y}(-Z)$ denotes the ideal sheaf of $Z$ in $Y$. Thus,
viewed as an element of $K_{T}(X)$,
$[\co_{Y}(-Z)]=[\co_{Y}]-[\co_{Z}]$. Moreover, if $\cf$ is a
$T$-equivariant coherent sheaf on $Y$ then $\cf(-Z)$ shall denote $\cf
\otimes \co_{Y}(-Z)$.

Recall that Demazure has defined the operators $L_w$, $w\in W$,
on
$\bz[\Lambda]$ satisfying the following properties:
\be\label{propdem}\begin{array}{llll} &L_wL_{w^{'}}&=L_{ww^{'}}
  &\mbox{if}~~l(ww^{'})=l(w)+l(w^{'})\\ &L_sL_s&=L_s &\mbox{if}
~~l(s)=1~;~i.e~ s=s_{\alpha}~ \mbox{for~some}~\alpha\in\Delta
\\&L_{s_{\alpha}}(f)&=\frac{f-{s_{\alpha}(f)}}{1-e^{\alpha}}&\mbox{for}~~\alpha\in
  \Delta~~\mbox{and}~~f\in R(T)\end{array}\ee (see Theorem 2,
pp. 86-87 of \cite{d}). In particular, $L_{s_{\alpha}}$ (and
hence
$L_w$) is a $R(T)^{W}$-linear operator on $R(T)$. Also, for a
simple
reflection $s\in W$ we have: \be\label{propdem2} L_s\cdot
L_w=\left\{\begin{array}{ll} L_{sw} & \mbox{if
$l(sw)=l(w)+1$}\\L_w &
    \mbox{if $l(sw)=l(w)-1$}
                        \end{array}
\right.\ee Moreover, for any $w'\in W$, there exists a unique
$v(w,w')\in W$ such that:
\be\label{bar{w}}L_{w^{\prime}}L_{w^{-1}w_0}=L_{v(w,w')}~\ee
(see \S5.6 of \cite{d}).

                      For each proper $T$-variety $Y$ and each
$T$-equivariant coherent sheaf $\cf$ on $Y$, we define
\be\label{eqeulerchar}
                      \chi^{T}(Y,[\cf])=\pi_{*}([\cf])\ee where
                      $\pi:Y\ra pt$ is the unique map to a
                      point. Observe that,
\be\label{ep}\chi^{T}(Y,[\cf])=\sum_{k}(-1)^{k}\mbox{Char~}
(H^{k}(Y,\cf ))\in R(T) \ee where
$\mbox{Char~}(H^{k}(Y,\cf))\in R(T)$ is the character of the
finite dimensional $T$-module
\[H^{k}(Y,\cf)=H^{k}(G/B,\co_{Y}\otimes
\cf).\] Further, $\chi^T:K_T(Y)\ra R(T)$ is an $R(T)$-linear map.
In particular, let $\cf$ be a $T$-equivariant coherent sheaf on
$G/B$. Then we define the equivariant Euler-Poincar\'{e}
characteristic
as:\be\label{epdem}\chi^{T}(X_w,c_{K}^{T}(e^{\lambda}))=e^{\rho}\cdot
                      L_{w} (e^{\lambda-\rho}).\ee Moreover, if
$\epsilon:R(T)\ra \bz$ denotes the canonical augmentation, then
\be\label{demchar}\chi(X_w,\cl(\lambda))=\epsilon L_w(e^{\lambda-\rho})\ee
where $\chi(.,.)$ denotes the ordinary Euler-Poincar\'{e}
characteristic (see Theorem 2(b) and Cor.1 pp. 86-87 of
\cite{d}).

                      In \cite{kk}, Kostant and Kumar define an
$R(T)$-module basis $(\tau^{w})_{w\in W}$ for $K_{T}(G/B)$ which satisfies:
\be\label{dual}\chi^{T}(X_{v^{-1}},*\tau^{w})=\delta_{v,w}\ee
(see p. 591, Prop.3.39 of \cite{kk}). Let
                      \be\label{boundoppsch} \partial
X^w:=\bigsqcup_{v\in W,~ v>w} B^{-}vB/B \ee
and\be\label{gkdualbasis}\xi^{w}:=[\co_{X^w}(-\partial
X^w)]_{T}.\ee Recall from Prop. 2.1 of
\cite{gk}, that $\{\xi^{w}\}_{w\in W}$ form an $R(T)$-basis for
$K_T(G/B)$ dual to the Schubert basis $\{[\co_{X_w}]_{T}\}_{w\in
W}$ under the pairing: \be\label{pairing}\langle u,
v\rangle:=\chi^{T}(G/B , u\cdot v) ~~~\mbox{forall ~$u,v\in
K_{T}(G/B)$}.\ee Further, it is shown in Prop. 2.2 of \cite{gk}
that:
\be\label{gk} *\tau^w=\xi^{w^{-1}}\ee where $\tau^w$ is the
Kostant-Kumar basis.

We further have the following relation between the Graham-Kumar
basis
and the opposite Schubert basis of $K_T(G/B)$ (see \cite{gk}):
\be\label{kks}[\co_{X^w}]_{T}=\sum_{w\preceq w^{\prime}}
\xi^{w^{\prime}}.\ee

For $I\subseteq \Delta$, let $W_I$ be the subgroup of $W$
generated by
$\{s_\alpha:\alpha\in I\}$. Further, let $W^I$ denote the
minimal
length coset representatives of $W/W_I$. Let $P=P_{I}\supset B$
denote
the corresponding standard parabolic subgroup. In particular,
for
$w\in W^{I}$, we have the Schubert variety $X_w^P$ (resp. the
opposite
Schubert variety $X^w_{P}$) which is the closure of the Bruhat
cell
$BwP/P$ (resp. opposite Bruhat cell $B^{-}wP/P$) in the partial
flag
variety $G/P$.

It is well known that $\{[\co_{X_w^{P}}]_{T}\}_{w\in W^{I}}$ is
an
$R(T)$-basis of $K_{T}(G/P)$, and so is
$\{[\co_{X^w_{P}}]_{T}\}_{w\in
W^I}$. Further, in \S2 of \cite{gk}, Graham and Kumar define
the elements: \be\label{gkdef}\xi_P^{v}=[\co_{X_{P}^v}(-\partial
X^v_{P})]_{T}\ee which form an $R(T)$-basis
$\{\xi_{P}^{v}\}_{v\in
  W^{I}}$ for $K_{T}(G/P)$ dual to the Schubert basis
$\{[\co_{X^{P}_w}]_{T}\}_{w\in W^{I}}$ under the pairing:
\be\label{pairingpartial}\langle [\co_{X^{P}_w}]_{T},
\xi_{P}^{v}\rangle =\chi^{T}(G/P,\co_{X^{P}_w\cap
X_{P}^v}(-X^{P}_w
\cap \partial X_{P}^v)).\ee Note that (\ref{pairingpartial}) is
a
generalization of (\ref{pairing}) to $G/P$.

Let $K(G/B)$ denote the Grothendieck ring of coherent sheaves
on $G/B$. Further,
\[c_{K}:\bz[\Lambda]=R(T)=K_{G}(G/B)\ra K(G/B)\] denote the
characteristic homomorphism. Recall that in \cite{d}, Demazure
has
established the existence of a basis $(a_{w})_{w\in W}$ for the
$\bz$-module $K(G/B)$ such that: \be\label{char}
c_{K}(e^{\lambda})=\sum_{w\in W} \chi(X_w,\cl(\lambda))a_w \ee
where
$\chi(.,.)$ denotes the Euler-Poincar\'{e} characteristic.

Recall that we have the following relation between the Demazure
basis
and Schubert basis (see Prop. 4.3.2 of \cite{br4}):
\be\label{ds}[\co_{X^w}]=\sum_{w\preceq
w^{\prime}} a_{w^{\prime}}.\ee

Let \be\label{forget} {f}:K_{T}(G/B)\ra K(G/B)\ee denote the
forgetful
homomorphism. Then we have \[{
f}([\co_{X^{w}}]_T)=[\co_{X^{w}}])\] and
${f}(\xi^{w})=a_w$ (see Prop. 3.39 of \cite{kk}).

{\bf Acknowledgement:} I am grateful to Prof. Michel Brion for
several valuable discussions and suggestions during this work. 
I also thank him for a careful reading and
invaluable comments on many earlier versions of this
manuscript.I
thank Prof. Shrawan Kumar for some motivating questions and for
sending me his paper with W. Graham (\cite{gk}) which was key
to this
work.

\section{Equivariant $K$-theory of flag varieties}

\subsection{Lifting of Schubert basis to $R(T)\otimes R(T)$}
In this section we construct explicit lifts of classes of
structure
sheaves of Schubert varieties in $K_{T}(G/B)$ to the ring
$R(T)\otimes
R(T)$. These shall be used in the later sections to describe
the multiplicative structure of $K_{T}(G/B)$.

{\it We mention here that tensor products are considered over
$\bz$
  unless specified.}

\blem\label{genreln} The canonical homomorphism
$$\Psi:R(T)\otimes R(T)\ra K_T(G/B)$$ that sends an element
$\sum_{i=1}^{n}a_i\otimes b_i$ in $R(T)\otimes R(T)$ to the
element\\ $\sum_{i=1}^{n}a_i\cdot c^{T}_{K}(b_i)$ in
$K_{T}(G/B)$ is
surjective with kernel the ideal $$\ci=\langle c\otimes
1-1\otimes
c~:~ c\in R(T)^{W}\rangle$$ in $R(T)\otimes R(T)$. \elem {\bf
Proof:}
We recall from Prop. 4.1 of \cite{mer} that the map:
\be\label{isomerk} R(T)\otimes_{R(G)}
K_{G}(G/B)=R(T)\otimes_{R(T)^W}
R(T)\rightarrow K_{T}(G/B),\ee defined via $a\otimes b~\mapsto
a\cdot
c^{T}_{K}(b)$ is an isomorphism.  Moreover, by definition of
$R(T)\otimes_{R(T)^W} R(T)$, there is canonical surjective
homomorphism $\psi: R(T)\otimes R(T)\ra
R(T)\otimes_{R(T)^{W}}R(T)$
with kernel precisely $\ci$. Now, $\Psi$ is the homomorphism
obtained
by composing $\psi$ with the isomorphism given by
(\ref{isomerk}). Thus it follows that $\Psi$ is surjective with
kernel
$\ci$.  (Also see Theorem 1.2 of \cite{gr}). $\Box$

\bdefe\label{extdemop} By defining $\mathbb{L}_{w}(a\otimes
b):=a\otimes L_w(b)$ and extending it by linearity, we can
define the
Demazure operator $\mathbb{L}_w$ on $R(T)\otimes R(T)$ as an
$R(T)\otimes 1$-linear operator.\edefe

We now prove a preliminary lemma which shall be applied in the
main proposition.

\blem\label{sublem1} Let $v(w,w')$ be as in (\ref{bar{w}}).
Then\be\label{imply}v(w,w')=w_0 \Leftrightarrow w\preceq
w^{\prime}.\ee
\elem {\bf Proof:} Let $l(w)=r$ and $w_r:=w^{-1}w_0$. Let
$w'=s'_1\cdots s'_k$ be a reduced expression for $w'$. Hence
\be\label{bar{w}1} L_{w'}\cdot L_{w^{-1}w_0}=L_{w'}\cdot
L_{w_r}=L_{s'_1}\cdots L_{s'_k}\cdot L_{w_r}.\ee Now, by
(\ref{propdem2}) we see that: \be\label{express} L_{s'_1}\cdots
L_{s'_k}\cdot L_{w_r}=L_{s'_{i_1}s'_{i_2}\cdots s'_{i_m}w_r}\ee
for
some subsequence $(s'_{i_1},\ldots,s'_{i_m})$ of
$(s'_1,\ldots,s'_k)$.
Now, (\ref{bar{w}}), (\ref{bar{w}1}) and (\ref{express}) imply
that
\be\label{barw2}v(w,w')=s'_{i_1}s'_{i_2}\cdots s'_{i_m}w_r.\ee
Since
$w_0=w\cdot w_r$, it follows that $v(w,w')=w_0$ will imply: \be
w=s'_{i_1}s'_{i_2}\cdots s'_{i_m}.\ee Hence $w\preceq w'$ (see
Cor. 2.2.2 of \cite{br4}).

For the converse, we need to show that: 

{\it Claim}: If $w\preceq w'$ then
$L_{w'}L_{w^{-1}w_0}=L_{w_0}$.  

{\it Proof of Claim:} Note that when $l(w)=0$ then $w=1$. Thus
by (\ref{propdem2}) we
have: \be L_{w'}L_{w^{-1}w_0}=L_{w'}L_{w_0}=L_{w_0}.\ee Also,
when
$l(w')-l(w)=0$ then $w=w'$. Again by (\ref{propdem2}) we have:
\be
L_{w'}L_{w^{-1}w_0}=L_wL_{w_r}=L_{w_0}.\ee We shall now prove
the
claim by induction on $l(w)$ and $l(w')-l(w)$.

Let $s'_1\cdots s'_k$ be a reduced expression of $w'$ and let
$w=s'_{i_1}\cdots s'_{i_m}$. Now, we can write:\be
\label{caseexp}
L_{w'}L_{w^{-1}w_0}=L_{v}L_{s'_k}L_{w^{-1}w_0}\ee where
$v=s'_1\cdots
s'_{k-1}$. 

\noindent
{\it Case (i)} If $l(s'_kw_r)=l(w_r)-1$ then by
(\ref{propdem2})we
have: \be\label{caseeq1}
L_{v}L_{s'_k}L_{w^{-1}w_0}=L_vL_{w_r}=L_vL_{w^{-1}w_0}.\ee
Moreover,
we note that $l(s'_kw_r)=l(w_r)-1$ is equivalent to
$l(ws'_k)=l(w)+1$.
This further implies that $i_m\leq k-1$. Hence it follows that
$w\preceq v$. Now, since $l(v)-l(w)\lneq l(w')-l(w)$, the claim
follows by
induction on $l(w')-l(w)$.

\noindent
{\it Case (ii)} If $l(s'_kw_r)=l(w_r)+1$ then again by
(\ref{propdem2}) we have:
\be\label{caseeq2}L_{v}L_{s'_k}L_{w^{-1}w_0}=L_vL_{s'_kw_r}=L_vL_{(ws'_k)
^{-1}w_0} .\ee Note that $l(s'_kw_r)=l(w_r)+1$ implies
$l(ws'_k)=l(w)-1$. Since $w\preceq w'$, this further implies
that
$ws'_k\preceq w'$ (see Proposition on p.119 of \cite{h}).
Moreover,
since $l(ws'_k)\lneq l(w)$, we further see that $ws'_k\preceq
v$.
The claim now follows by induction on
$l(w)$. 

This completes the proof.$\Box$

\bpropo\label{slift} In $R(T)\otimes R(T)$ there exists an element
$u_{0}$ such that \be\label{sl}
\Psi(\mathbb{L}_{w^{-1}w_{0}}(u_{0})\cdot (1\otimes
e^{\rho}))=[\co_{X^w}]_{T}.\ee Indeed, we may take $u_0=v_0(1\otimes
e^{-\rho})$ where $v_0$ is such that $\Psi(v_0)=[\co_{X^{w_0}}]_{T}$.
\epropo {\bf Proof:} By (\ref{dual}) and (\ref{gk}), we have the
following identity in $K_{T}(G/B)$: \be\label{eqchar}
c^{T}_{K}(e^{\lambda})=\sum_{w\in W} \chi^{T}(X_w,\cl(\lambda))\xi^w
.\ee Moreover, combining (\ref{epdem}) and (\ref{eqchar}) it also
follows that: \be\label{eqchardem} c^{T}_{K}(e^{\lambda})=\sum_{w\in
  W}e^{\rho}\cdot L_{w} (e^{\lambda-\rho}) \xi^w .\ee By (\ref{kks}),
it follows in particular that
$\xi^{w_0}=[\co_{X^{w_0}}]_{T}=w_0\cdot[\co_{X_1}]_{T}$.

Now, since $\Psi$ is surjective by Lemma \ref{genreln}, there
exists
an element $v_0$ such that $\Psi(v_0)=\xi^{w_0}$. More
precisely, if
\be\label{not} v_0=\sum_{i=1}^{n}a_i\otimes b_i,\ee then using
(\ref{eqchardem}) we have:
\be\label{eq1}\Psi(v_0)=\sum_{i=1}^{n}a_i\cdot
c_{K}^{T}(b_i)=\sum_{i=1}^n\cdot\sum_{w\in W}a_i\cdot
e^{\rho}\cdot
L_{w}(b_i\cdot e^{-\rho})\xi^{w}.\ee Hence we see that:
\be\label{eq2} \Psi(v_0)=\sum_{w\in W}\sum_{i=1}^{n}a_i\cdot
e^{\rho}\cdot L_{w}(b_i\cdot e^{-\rho})\xi^{w} =\xi^{w_0}.\ee
Therefore we have: \be\label{eq3}\sum_{i=1}^{n}a_i\cdot
e^{\rho}\cdot
L_{w}(b_i\cdot e^{-\rho})=\delta_{w,w_0}.\ee

Let \be\label{notn} u_0:=v_0\cdot (1\otimes
e^{-\rho}).\ee 

{\it Claim:} $u_0$ is the required element in $R(T)\otimes
R(T)$that
satisfies (\ref{sl}).

{\it Proof of Claim:} Note that if $v_0$ is as in (\ref{not})
then:
\be\label{notn1}u_0=\sum_{i=1}^{n}a_i\otimes e^{-\rho}\cdot b_i
.\ee

Now, by (\ref{notn1}) and Def. \ref{extdemop} it
follows that: \be\label{reexp1}
\mathbb{L}_{w^{-1}w_0}(u_0)\cdot(1\otimes e^{\rho})=\sum_{i=1}^na_i\otimes e^{\rho}\cdot
L_{w^{-1}w_0}(e^{-\rho}\cdot b_i).\ee

Hence by (\ref{eqchardem}) and (\ref{reexp1}) we have:
\be\label{exp}
\Psi(\mathbb{L}_{w^{-1}w_0}(u_0)\cdot (1\otimes
e^{\rho}))=\sum_{w^{\prime}\in W} \sum_{i=1}^ne^{\rho}\cdot
a_i\cdot
{L}_{w^{\prime}}{L}_{w^{-1}w_0}( b_i\cdot
e^{-\rho})\xi^{w^{\prime}}.\ee Now we see that the claim
follows by (\ref{kks}),
(\ref{eq3}), (\ref{exp}) and Lemma \ref{sublem1}.$\Box$

\blem\label{slift2pr} If $w\in W^I$ and $r\in R(T)$, then we
have:
\be e^{\rho}\cdot L_{w^{-1}w_0}(e^{-\rho}\cdot r)\in
R(T)^{W_I}.\ee
\elem

\noindent
{\bf Proof:} We first note that:\be
\label{step2} s_j(e^{\rho} \cdot L_{w^{-1}w_{0}}(e^{-\rho}\cdot
r))=e^{\rho-\alpha_j}\cdot
s_j(L_{w^{-1}w_{0}}(e^{-\rho}\cdot r)).\ee

Thus we see that for $j\in I$, the condition:
\be\label{step1}s_j(e^{\rho}\cdot
L_{w^{-1}w_{0}}(e^{-\rho}\cdot r))=e^{\rho}\cdot
L_{w^{-1}w_{0}}(e^{-\rho}\cdot r)\ee is
equivalent
to: \be
\label{step3} s_j(L_{w^{-1}w_{0}}(e^{-\rho}\cdot
r))=e^{\alpha_j}\cdot L_{w^{-1}w_{0}}(e^{-\rho}\cdot r).\ee
Further, note that \be\label{step4}
s_j(L_{w^{-1}w_{0}}(e^{-\rho}\cdot
r))= L_{w^{-1}w_{0}}(e^{-\rho}\cdot r)-(1-e^{\alpha_j})\cdot
L_{s_j}L_{w^{-1}w_{0}}(e^{-\rho}\cdot r).\ee Let
$w_1:=w^{-1}w_0$. Then we have \be\label{step5}
l(s_jw_1)=l(s_jw^{-1}w_0)=l(w_0)-l(s_jw^{-1}).\ee Now, if $w\in
W^I$,
then for every $j\in I$ we have $l(ws_j)=l(w)+1$ which is
equivalent
to $l(s_jw^{-1})=l(w^{-1})+1$. This implies by (\ref{step5})
that:
\be\label{step6}\begin{array}{ll}l(s_jw_1)&=
l(w_0)-(l(w^{-1})+1)
  \\ &=l(w_0)-l(w^{-1})-1 \\&=l(w^{-1}w_0)-1\\&=l(w_1)-1.
\end{array}\ee This further implies by (\ref{propdem2}) 
that:
\be\label{step7}L_{s_j}L_{w^{-1}w_0}=L_{s_j}L_{w_1}=L_{w_1}.\ee

Now by substituting (\ref{step7}) in (\ref{step4}), we see that
when
$w\in W^I$, the condition (\ref{step3}) and hence (\ref{step1})
hold
for all $j\in I$. This proves that if $w\in W^I$ then
$e^{\rho}\cdot
L_{w^{-1}w_{0}}(e^{-\rho}\cdot r)\in R(T)^{W_I}$.$\Box$.

\bpropo\label{slift2} Let $u_0$ be as in Prop. \ref{slift}.
Then the
element $$\mathbb{L}_{w^{-1}w_{0}}(u_{0})\cdot (1\otimes
e^{\rho})$$
belongs to $R(T)\otimes R(T)^{W_I}$ if $w\in W^I$. \epropo 

{\bf Proof :} This proposition follows immediately from (\ref{reexp1})
and Lemma \ref{slift2pr}.$\Box$.
 
\bnot In the following sections we let $u_0=\sum_{i=1}^n a_i\otimes
e^{-\rho}\cdot b_i\in R(T)\otimes R(T)$ be as in
Prop.\ref{slift}. (See \S5 about this choice.)\enot

\subsection{Structure constants of Schubert basis in
$K_{T}(G/B)$}
In this section we determine a closed formula for the multiplicative
structure constants of the basis $\{[\co_{X^{w}}]_{T}\}_{w\in W}$ in
$K_{T}(G/B)$ in terms of the above elements $a_i, b_i$. We remark here
that in \cite{gk}, these structure constants as well as those of the
dual basis have been studied in detail with regard to the positivity
conjectures viz., Conjecture 3.1 and 3.10 of \cite{gk}. The author is
currently working to find more direct inter-connections between the
results in this section and those in \cite{gk}.

\blem\label{Multschubert} For $x,y,z\in W$, let
\be\label{strconst1eq}
C^z_{x,y}:= \sum_{w\preceq z}(-1)^{l(z)-l(w)} \sum_{1\leq
i,j\leq n}
a_i\cdot a_j\cdot e^{\rho}\cdot L_w(L_{x^{-1}w_{0}}(b_{i}\cdot
e^{-\rho})\cdot L_{y^{-1}w_{0}}(b_j\cdot e^{-\rho})\cdot
e^{\rho}) \ee
where $v_0=\sum^n_{i=1} a_i\otimes b_i\in R(T)\otimes R(T)$ is
such
that $\Psi(v_0)=[\co_{X^{w_0}}]_{T}$. Then in $K_{T}(G/B)$ we
have:
\be\label{lreq}[\co_{X^x}]_{T}[\co_{X^y}]_{T}=\sum_{z\in
W}C^z_{x,y}[{\mathcal O}_{X^{z}}]_{T}\ee for $x,y\in W$.  \elem
{\bf Proof:} Recall from Lemma 4.2 of \cite{gk} that the basis
$\{\xi^{v}\}_{v\in W}$ can be expressed in terms of the
Schubert basis
$\{[\co_{X^v}]_{T}\}_{v\in W}$ in $K_T(G/B)$ as follows:
\be\label{ksb}\xi^v=\sum_{v\preceq w}
(-1)^{l(w)-l(v)}[\co_{X^w}]_{T}.\ee Note that (\ref{ksb}) is
equivalent to (\ref{kks}) via M$\ddot{\mbox{o}}$bius inversion
(see Remark 4.3.3 of \cite{br4}). Now,
using Lemma \ref{genreln} and substituting (\ref{ksb}) in
(\ref{eqchardem}) we get: \be\label{schareq} \Psi(a\otimes
b)=\sum_{w\in W}\sum_{v\in W,v\preceq w}(-1)^{l(w)-l(v)} a\cdot
e^{\rho}\cdot
{L}_{v}(b\cdot e^{-\rho})[\co_{X^{w}}]_{T} \ee for $a\otimes
b\in
R(T)\otimes R(T)$.  

Moreover, by (\ref{reexp1}),
$\mathbb{L}_{x^{-1}w_{0}}(u_{0})\cdot(1\otimes
e^{\rho})\cdot\mathbb{L}_{y^{-1}w_{0}}(u_{0})\cdot(1\otimes
e^{\rho})=$\be\label{prodexp}\sum_{1\leq i,j\leq n}a_i\cdot
a_j\otimes
e^{2\rho}\cdot L_{x^{-1}w_0}(e^{-\rho}\cdot b_i)\cdot
L_{y^{-1}w_0}(e^{-\rho}\cdot b_j).\ee

Further, by Prop. \ref{slift} it follows that:
\be\label{sprodeq}\Psi(\mathbb{L}_{x^{-1}w_{0}}(u_{0})\cdot
(1\otimes
e^{\rho})\cdot \mathbb{L}_{y^{-1}w_{0}}(u_{0})\cdot (1\otimes
e^{\rho}))=[\co_{X^x}]_{T}[\co_{X^y}]_{T}\ee for $x,y\in W$.
Then by
(\ref{sprodeq}) and (\ref{schareq}) we get (\ref{lreq}) where
$C^{z}_{x,y}$ is as in (\ref{strconst1eq}). $\Box$

\subsection{A Chevalley formula in $K_{T}(G/B)$}
The following lemma gives ``a Chevalley formula'' in
$K_{T}(G/B)$,
which determines the coefficients when the product
$[\cl^{T}(\lambda)]_{T}[\co_{X^{x}}]_{T}$ is expressed in terms
of the
Schubert basis $\{[\co_{X^v}]_{T}:{v\in W}\}$.

\blem\label{chevalley} For $\lambda\in X^{*}(T)$ and $x,y\in W$
let
\be\label{lbscheq} Q^{\lambda}_{x,y}:=\sum_{w\in W,w\preceq
  y}(-1)^{l(y)-l(w)}\sum_{i=1}^n e^{\rho}\cdot a_i \cdot
L_w(e^{\lambda}\cdot L_{x^{-1}w_0}(b_i\cdot e^{-\rho}))\ee
where $v_0=\sum^n_{i=1} a_i\otimes b_i\in R(T)\otimes R(T)$ is
such
that
$\Psi(v_0)=[\co_{X^{w_0}}]_{T}$. Then in $K_{T}(G/B)$ we have:
\be\label{cheveq}[\cl^{T}(\lambda)]_{T}\cdot
[\co_{X^{x}}]_{T}=\sum_{y\in
W}Q^{\lambda}_{x,y}[\co_{X^y}]_{T}.\ee
\elem {\bf Proof:} Note that \be\label{liftlb}\Psi(1\otimes
e^{\lambda})=c^{T}_{K}(e^{\lambda})=[\cl^{T}(\lambda)]_{T}.\ee
By
(\ref{liftlb}) and Prop.\ref{slift} it follows that:
\be\label{prodlbsch}\Psi(\mathbb{L}_{x^{-1}w_{0}}(u_{0})\cdot
(1\otimes e^{\rho})\cdot (1\otimes
e^{\lambda}))=[\co_{X^x}]_{T}\cdot
[\cl^{T}(\lambda)]_{T}\ee for $x\in W$ and $\lambda\in
X^{*}(T)$. Thus
we see that (\ref{cheveq}) follows immediately from
(\ref{reexp1}),
(\ref{schareq}) and (\ref{prodlbsch}) where $Q^{\lambda}_{x,y}$
is
given by (\ref{lbscheq}). $\Box$

\subsubsection{Comparison with known Chevalley formulas}

\blem\label{cohomchar} Let $x,y\in W$ and $w=w_0x$, $v=w_0y$.
We then
have the following interpretation of (\ref{lbscheq}):
\be\label{lbscheq1}Q^{\lambda}_{x,y}= w_0\cdot
\chi^{T}(X^{v}_{w},\cl^{T}(w_0(\lambda))(-(\partial
X^{v})_w))\ee
whenever $x\preceq y$ and $Q^{\lambda}_{x,y}=0$ otherwise (see
Lemma 1 of \cite{br3}). In
particular, when $w_0(\lambda)\in X^{*}(T)$ is dominant we
have:\be\label{lbscheq2} Q^{\lambda}_{x,y}=w_0\cdot
\mbox{Char~}H^0(X^{v}_{w},\cl^{T}(w_0(\lambda))(-(\partial
X^{v})_w))\ee whenever $x\preceq y$ and $Q^{\lambda}_{x,y}=0$
otherwise.\elem {\bf Proof:} Let $\xi_{y}:=[\co_{X_y}(-\partial
X_y)]$. Then $\xi_y$ is dual to $[\co_{X^x}]$ under the pairing
(\ref{pairing}).
Now,  (\ref{cheveq}) implies that:
\be\label{lbsch2eq}\begin{array}{ll}Q^{\lambda}_{x,y}&=\langle
  ~[\cl^{T}(\lambda)]_{T}\cdot [\co_{X^x}]_{T}~,~
  \xi_{y}~\rangle\\&=\chi^{T}(X^{x}_{y},~
[\cl^{T}(\lambda)]_{T}(-\partial X_{y})^{x}]\end{array}\ee
whenever
$x\preceq y$, and $Q^{\lambda}_{x,y}=0$ otherwise. The second
equality above follows because the intersections $X_y\cap X^x$
and $X^x\cap\partial{X_y}$ are transversal (see Lemma 4.1.2 of \cite{br4}).

If $w=w_0x$ and $v=w_0y$, we can write (\ref{lbsch2eq}) as:
\be\label{lbsch3eq}\begin{array}{ll}Q^{\lambda}_{x,y}&=\chi^{T}(X^{w_0w}_{w_0v},\cl^{T}(\lambda)(-(\partial
  X_{w_0v})^{w_0w}))\\&=\chi^{T}(w_0\cdot
  X^{v}_{w},\cl^{T}(\lambda)(-w_0\cdot(\partial
  X^{v})_w))\\&=w_0\cdot\chi^{T}(
  X^{v}_{w},\cl^{T}(w_0(\lambda))(-(\partial X^{v})_w)).
\end{array}\ee whenever $v\preceq w$ and$Q^{\lambda}_{x,y}=0$
otherwise. (Note that $x\preceq y$ is
equivalent
to $v\preceq w$.) In particular, when $w_0(\lambda)\in X^*(T)$
is
dominant, by Prop.1 on p.9 of \cite{bl} it follows that:
$$\chi^{T}(X^{v}_{w},\cl^{T}(w_0(\lambda))(-(\partial
X^{v})_w))=\mbox{Char~}H^0(X^{v}_{w},\cl^{T}(w_0(\lambda))(-(\partial
X^{v})_w)).$$ Hence the lemma.$\Box$

\brem\label{compare} Let $w=w_0x$ and $v=w_0y$. Then
(\ref{cheveq})
can be rewritten as: \be\label{chev1eq}[\cl^{T}(w_0
\lambda)]_{T}\cdot
[\co_{X_{w}}]_{T}=\sum_{v\preceq w}w_0(
Q^{\lambda}_{x,y})[\co_{X_{v}}]_{T}\ee where
$Q^{\lambda}_{x,y}$is as
in (\ref{lbsch3eq}).  Hence substituting (\ref{lbsch3eq}) in
(\ref{chev1eq}), we derive the ``Chevalley formula'' as in
\cite{ls}
and \cite{ll}. In particular, note that
$w_0(Q^{w_0\lambda}_{w_0w,w_0v})$ is same as
$C^{\lambda}_{w,v}$of
\cite{ls} where it is interpreted as $\sum e^{-\pi(1)}$ where
the sum
runs over all L-S paths $\pi$ of shape $\lambda$ ending in $v$
and
starting with an element smaller or equal to $w$. Here we
briefly
recall that an $L-S$ path $\pi$ of shape $\lambda$ on $X(\tau)$
is a
pair of sequences $\pi=(\underline{\tau},\underline{a})$ of
Weyl group
elements and rational numbers, where $\underline{\tau}$ is of
the form
$\underline{\tau}=(\tau_1,\ldots,\tau_r)$ such that $\tau\geq
\tau_1$
and $\tau_1\geq \tau_2\geq\cdots\geq\tau_r$ in the Bruhat order
on
$W$. We call $\tau_1=i(\pi)$ the ``initial element of $\pi$ and
$\tau_r=e(\pi)$ the ``end'' element of $\pi$. (see \S3 of
\cite{ls}). \erem

\brem\label{otherref} We refer to \cite{ls} and \cite{ll} for
more
details on representation theoretic interpretation of the
Chevalley
formula using Standard Monomial Theory. We also refer to
\cite{gr} for
Chevalley formula in $K_{T}(G/B)$ given in terms of the
combinatorics
of the Littelmann path model, using affine nil
Hecke-Algebra. Also see \cite{mw} for recent results on
Chevalley
formula in equivariant $K$-theory of flag varieties using the
Bott-Samelson resolution.  \erem

\subsection{Structure constants of Schubert basis in
$K_{T}(G/P)$}

In this section we determine a closed formula for the
multiplicative
structure constants of the Schubert basis
$\{[\co_{X^{w}_{P}}]_{T}\}_{w\in W^{I}}$ of $K_{T}(G/P)$, again
in
terms of $a_i,b_i$.

Let $\mu^{I}$ be the M$\ddot{\mbox{o}}$bius function of the
induced
Bruhat ordering on $W^I$. Then (see Theorem 1.2 of \cite{de}):
\be\label{relmob} \mu^I(v,w)=\left\{\begin{array}{ll}
(-1)^{l(v)+l(w)}
& \mbox{if $[v,w]\cap W^I=[v,w]$}\\ 0 &
\mbox{otherwise}\end{array}\right.\ee where for $v\preceq w$,
$[v,w]:=\{u\in W : v\preceq u\preceq w\}$.

\blem\label{Multpartial} For $x,y,z\in W^I$, let
\be\label{strconst2eq} D^z_{x,y}:=
 \sum_{w\in W^I, w\preceq z}\mu^{I}(w,z)
\sum_{1\leq i,j\leq n} e^{\rho}\cdot a_i\cdot a_j \cdot
L_w(L_{x^{-1}w_{0}}(b_i\cdot e^{-\rho})\cdot
L_{y^{-1}w_{0}}(b_j\cdot
e^{-\rho})\cdot e^{\rho}). \ee Then in $K_{T}(G/P)$ we have:
\be\label{lr2eq}[\co_{X^x_P}]_{T}[\co_{X^y_P}]_{T}=\sum_{z\in
W^I}D^z_{x,y}[{\mathcal O}_{X^{z}_P}]_{T}\ee for $x,y\in W^I$.
\elem

{\bf Proof:} Let $\pi:G/B \ra G/P$ be the canonical projection.
Then
we have: \be\label{liftshu}
\pi^{*}([\co_{X^w_{P}}]_{T})=[\co_{X^w}]_{T}~~~ \mbox{for~
$w\in W^{I}$,}
\ee where $\pi^{*}: K_{T}(G/P)\ra K_{T}(G/B)$ is the induced
morphism.

For any $v\in W^{I}$ we also have (see Lemma 3.4 of \cite{gk}):
\be\label{sumgk}\pi^{*}\xi_{P}^v=\sum_{u\in W_I}\xi^{vu}.\ee
Furthermore, for $w,v\in W^{I}$, we shall identify the elements
$\xi^v_{P}$ and $[\co^w_{P}]_{T}$ in $K_{T}(G/P)$ with their
images in
$K_{T}(G/B)$ under the injective morphism $\pi^*$.

Further, it follows from (\ref{propdem}) that for $r\in
R(T)^{W_{I}}$
and $\alpha\in I$ we have: \be\label{marlindem}
L_{s_{\alpha}}(r\cdot
e^{-\rho})=\frac{r\cdot e^{-\rho}-r\cdot
e^{-\rho+\alpha}}{1-e^{\alpha}}=r\cdot e^{-\rho}.\ee This
implies
that for any $(w',v)\in W_I\times W^{I}$ we have:
\be\label{marlindem1}L_{vw'}(r\cdot e^{-\rho})=L_vL_{w'}(r\cdot
e^{-\rho})=L_v(r\cdot e^{-\rho}).\ee
Now, by (\ref{eqchardem}), (\ref{sumgk}) and (\ref{marlindem})
it
follows that for any $r\in R(T)^{W_I}$: \be\label{sdeq}
c^{T}_{K}(r)=\sum_{v\in W^I}e^{\rho} \cdot L_v(r\cdot
e^{-\rho})\cdot
\xi_{P}^v.\ee

Further, by Prop. \ref{slift} and (\ref{liftshu}) we have:
\be\label{sliftP} \Psi(\mathbb{L}_{w^{-1}w_{0}}(u_{0})\cdot
(1-e^{\rho}))=[\co_{X^w_P}]_{T}.\ee

By Lemma \ref{slift2pr} we have:
\[L_{w^{-1}w_{0}}(e^{-\rho}\cdot b_i)\cdot e^{\rho}\in
R(T)^{W_I}~~\mbox{for}
~~1\leq i\leq n\] whenever $w\in W^I$.

Now, from (\ref{exp}), (\ref{sumgk}) and substituting
$L_{w^{-1}w_{0}}(e^{-\rho})\cdot b_i\cdot e^{\rho}$ for $r$ in
(\ref{marlindem1}) we get: \be\label{sd1}[\co_{X^w_P}]_{T}=
\Psi(\mathbb{L}_{w^{-1}w_{0}}(u_{0})\cdot (1\otimes
e^{\rho}))=\sum_{v\in W^I}\sum_{i=1}^{n} e^{\rho}\cdot a_i\cdot
L_v(L_{w^{-1}w_{0}}(e^{-\rho}\cdot b_i))\xi^{v}_{P}.\ee Let
$v(w,w')\in W$ be such that $L_{v(w,w')}=L_v\cdot L_{w^{-1}w_{0}}$.
Then by (\ref{eq3}) we have: \be\label{recallprop}\sum_{i=1}^n
e^{\rho}\cdot a_i\cdot L_{v(w,w')}(e^{-\rho}\cdot
b_i)=\delta_{v(w,w'),w_0}.\ee Further, by Lemma \ref{sublem1},
(\ref{sd1}) can be rewritten as: \be\label{sfd} [\co_{X^w_P}]_{T}=
\Psi(\mathbb{L}_{w^{-1}w_{0}}(u_{0})\cdot (1\otimes
e^{\rho}))=\sum_{v\in W^{I}~,~w\preceq v} \xi^{v}_{P}.\ee

Further, via M$\ddot{\mbox{o}}$bius inversion (\ref{sfd}) is
equivalent to: \be\label{ksp}\xi_{P}^v=\sum_{w\in
W^{I}}\mu^I(v,w)[\co_{X_{P}^w}]_{T}\ee where $\mu^I(v,w)$ is
as defined in (\ref{relmob}).

Now, substituting (\ref{ksp}) in (\ref{sdeq}) and using Lemma
\ref{genreln} we get: \be\label{schareqpart} \Psi(t\otimes
r)=\sum_{w\in W^{I}} \mu^I(v,w)\cdot e^{\rho}\cdot t\cdot
{L}_{v}(r\cdot e^{-\rho})[\co_{X_{P}^{w}}]_{T} \ee for
$t\otimes r\in
R(T)\otimes R(T)^{W_{I}}$.  Since, by Prop. \ref{slift} and
Prop.\ref{slift2} we have:
$$[\co_{X^w_P}]=\Psi(\mathbb{L}_{w^{-1}w_{0}}(u_{0})\cdot
(1\otimes
e^{\rho})),$$ (\ref{schareqpart}) implies that:
\be\label{sprod1eq}\Psi(\mathbb{L}_{x^{-1}w_{0}}(u_{0})\cdot
(1\otimes
e^{\rho})\cdot \mathbb{L}_{y^{-1}w_{0}}(u_{0})\cdot (1\otimes
e^{\rho}))=[\co_{X^x_P}]_{T}[\co_{X^y_P}]_{T}\ee for $x,y\in
W^{I}$.

Thus by (\ref{prodexp}), (\ref{schareqpart}) and
(\ref{sprod1eq}) we
get (\ref{lr2eq}), where $D^{z}_{x,y}$ is as in
(\ref{strconst2eq}). Hence the lemma. $\Box$

\brem\label{generalise} Note that Lemma \ref{Multpartial} is a
generalization of Lemma \ref{Multschubert} to partial flag
varieties.
\erem

\section{Analogous results in ordinary $K$-theory}

In this section we construct explicit lifts of structure sheaves of
Schubert varieties in $K(G/B)$ in $1\otimes R(T)$. Indeed, the
forgetful homomorphism $f: R(T)\otimes_{R(T)^{W}} R(T)=K_{T}(G/B)\ra
K(G/B)$ lifts to a map $\wt{f}: R(T)\otimes R(T) \ra R(T)$ given by
$e^{\lambda}\otimes e^{\mu}\mapsto e^{\mu}$.

Let $v_0$ be as in (\ref{not}) and $u_0=(1\otimes e^{-\rho})\cdot
v_0$. Further, let $v'_0:=\wt{f}(v_0)$ and $u'_0:=\wt{f}(u_0)$.  Then
we have \be\label{imnot} v'_0=\sum_{i=}^n \epsilon(a_i)\cdot b_i .\ee
and \be\label{imnot1} u'_0=\sum_{i=1}^n \epsilon(a_i) \cdot
e^{-\rho}\cdot b_i.\ee

  The following proposition describes explicit lifts of
  $[\co_{X^w}]\in K(G/B)$ in $R(T)$.

  \bpropo\label{oslift} Let $v'_0, u'_{0}\in R(T)$ be as in
  (\ref{imnot}) and (\ref{imnot1}) respectively. Then
  $c_K(v'_0)=[\co_{X^{w_0}}]$ and \be\label{osl}
  c_{K}(L_{w^{-1}w_{0}}(u'_{0})\cdot e^{\rho})=[\co_{X^w}].\ee \epropo
  {\bf Proof:} Recall that $f([\co_{X^w}]_{T})= [\co_{X^w}]$ (see
  \ref{forget}). Moreover, by Prop. \ref{slift},
  $\mathbb{L}_{w^{-1}w_{0}}(u_{0})\cdot (1\otimes e^{\rho})$ lifts
  $[\co_{X^w}]_{T}$ in $R(T)\otimes R(T)$. Now, (\ref{reexp1}) implies
  that \[\wt{f}(\mathbb{L}_{w^{-1}w_{0}}(u_{0}))=
  L_{w^{-1}w_{0}}(u'_{0}).\] Thus by a simple diagram chase it follows
  that the image $L_{w^{-1}w_{0}}(u'_{0})\cdot e^{\rho}$ of
  $\mathbb{L}_{w^{-1}w_{0}}(u_{0})\cdot (1\otimes e^{\rho})$ under
  $\wt{f}$ lifts $[\co_{X^w}]$ in $R(T)$. $\Box$

\bpropo\label{oslift2} The element $L_{w^{-1}w_{0}}(u'_{0})\cdot
e^{\rho}\in R(T)^{W_I}$ if $w\in W^I$.  \epropo

\noindent
{\bf Proof :} This follows immediately by applying Lemma
\ref{slift2pr} for $r=v'_0$.$\Box$

We now state the results in ordinary $K$-ring analogous to Lemmas
\ref{Multschubert}, \ref{chevalley}, \ref{Multpartial}. Since the
proofs follow the same lines as the proofs of the above mentioned
lemmas, we avoid the repetition here.

\blem\label{oMultschubert} For $x,y,z\in W$, let
\be\label{ostrconst1} c^z_{x,y}:= \sum_{w\preceq
z}(-1)^{l(z)-l(w)}
\epsilon L_w(L_{x^{-1}w_{0}}(u'_0)\cdot
L_{y^{-1}w_{0}}(u'_0)\cdot e^{\rho}) .\ee Then in $K(G/B)$ we
have:
\be\label{olreq}[\co_{X^x}][\co_{X^y}]=\sum_{z\in
W}c^z_{x,y}[{\mathcal
    O}_{X^{z}}]\ee for $x,y\in W$.  \elem

\blem\label{ochevalley} For $\lambda\in X^{*}(T)$ and $x,y\in
W$let
\be\label{olbscheq} q^{\lambda}_{x,y}:=\sum_{w\preceq
  y}(-1)^{l(y)-l(w)}\epsilon L_w(e^{\lambda}\cdot
L_{x^{-1}w_0}(u'_0)).\ee Then in $K(G/B)$ we have:
\be\label{ochev}[\cl(\lambda)][\co_{X^{x}}]=\sum_{y\in
  W}q^{\lambda}_{x,y}[\co_{X^y}].\ee \elem

\blem\label{oMultpartial} For $x,y,z\in W^I$, let
\be\label{ostrconst2eq} d^z_{x,y}:= \sum_{w\preceq
z}\mu^{I}(z,w)
\epsilon L_w(L_{x^{-1}w_{0}}(u'_0)\cdot
L_{y^{-1}w_{0}}(u'_0)\cdot e^{\rho}). \ee Then in $K(G/P)$ we
have: \be\label{olr2eq}[\co_{X^x_P}][\co_{X^y_P}]=\sum_{z\in
W^I}d^z_{x,y}[{\mathcal O}_{X^{z}_P}]\ee for $x,y\in W^I$. \elem

\section{$K$-ring of the wonderful compactification}
\subsection{Some preliminaries}
Let $\alpha_1,\ldots,\alpha_r$ be an ordering of the set
$\Delta$ of
simple roots and $\omega_1,\ldots,\omega_r$ denote the
corresponding
fundamental weights for the root system of $(G,T)$. Since $G$
is simply connected, the fundamental weights form a basis for
$X^*(T)$
and hence for every $\lambda\in \Lambda$, $e^{\lambda}\in R(T)$
is a
Laurent monomial in the elements $e^{\omega_i}:1\leq i\leq r$.

In Theorem 2.2 of \cite{st} Steinberg has defined a basis
\be\label{steinbergbasis}\{f_v: v\in W^{I}\}\ee of
$R(T)^{W_{I}}$ as a
free $R(G)$-module of rank $|W^{I}|$. We recall here this
definition:
For $v\in W^{I}$ let
\[p_{v}:=\prod_{v^{-1}\alpha_{i}<0} e^{\omega_{i}}\in
R(T).\] Then
\[f_v:=\sum_{x\in W_{I}(v)\big{\backslash} W_{I}}
x^{-1}v^{-1}p_{v}\] where $W_{I}(v)$ denotes the stabilizer of
$v^{-1}p_v$ in $W_{I}$.  

Let $c_K:R(T)^{W_{I}}\ra K(G/P_I)$ denote the restriction of
the characteristic homomorphism (see \S8 of \cite{m}). Let
$I(G):=\{a-\epsilon(a)~|~a\in R(G)\}$ denote the augmentation
ideal. Then it is known that $c_{K}$ is surjective ring
homomorphism
and \be\label{aug} \mbox{ker}(c_K)=I(G)\cdot R(T)^{W_{I}}.\ee

Let $r_v:=L_{v^{-1}w_{0}}(u_{0})\cdot e^{\rho}\in R(T)^{W_I}$
for
${v\in W^I}$ for every $I\subset \Delta$. Then by Lemma
\ref{oslift2}
we recall that $c_{K}(r_v)=[\co_{X_{P_I}^v}]$ for every $v\in
W^{I}$.

Further, we recall that $R(T)^{W_{I}}=\bz[\Lambda]^{W_{I}}$ and
$R(G)=R(T)^{W}=\bz[\Lambda]^{W}$.

Note that $\mathfrak{p}:=I(G)$ is a prime ideal in $R(G)$ and
let
$R(G)_{\mathfrak{p}}$ denote the corresponding localization. We
further observe that the augmentation extends to
\be\label{extaug}
R(G)_{\mathfrak{p}}\ra \bq.\ee with kernel the maximal ideal
$\mathfrak{p}\cdot R(G)_{\mathfrak{p}}$. Further, the
characteristic
homomorphism extends to
\be\label{extchar}c_{K}:R(T)_{\mathfrak{p}}^{W_{I}}\twoheadrightarrow
K(G/P_{I})_{\bq}\ee with kernel $\mathfrak{p}\cdot
R(T)_{\mathfrak{p}}^{W_{I}}$.

\blem\label{basis2} The elements $\{ r_v:v\in W^{I}\}$ form a
basis of
$R(T)_{\mathfrak{p}}^{W_{I}}$ as an $R(G)_{\mathfrak{p}}$
module. \elem
           
{\bf Proof:} Now, $R(T)_{\mathfrak{p}}^{W_I}$ is a finitely
generated
$R(G)_{\mathfrak{p}}$-module. Moreover,
$\{c_{K}(r_v)=[\co_{X^v}]$, $v\in
W^{I}\}$ form a basis of $K(G/P_{I})_{\bq}$ as a
$R(G)_{\mathfrak{p}}/{\mathfrak{p}\cdot
R(G)_{\mathfrak{p}}}\simeq\bq$-vector space. Now, by Nakayama
lemma
(see Prop. 2.8 of \cite{am}) we see that $\{r_v:v\in W^{I}\}$
span
$R(T)_{\mathfrak{p}}^{W_I}$ as an
$R(G)_{\mathfrak{p}}$-module. Furthermore, since
$R(T)_{\mathfrak{p}}^{W_I}$ is free over
$R(G)_{\mathfrak{p}}$ of rank $|W^{I}|$, it follows that
$\{r_v:v\in W^{I}\}$ form a basis of
$R(T)_{\mathfrak{p}}^{W_I}$ as an
$R(G)_{\mathfrak{p}}$-module.$\Box$ 

We now fix some notations (also see p.378 of \cite{u}).

Note that $J\subseteq I$ implies that $W^{\Delta\setminus
J}\subseteq
W^{\Delta\setminus I}$. Let \be\label{1}
C^{I}:=W^{\Delta\setminus
I}\setminus (\bigcup_{J\subsetneq I}W^{\Delta \setminus J})\ee
and
\be\label{2} R(T)_{I}:=\bigoplus_{v\in
  C^I}R(T)_{\mathfrak{p}}^{W}\cdot r_v .\ee where
$R(T)_{\mathfrak{p}}^{W}=R(G)_{\mathfrak{p}}$.

\blem\label{dirsum} We have the following direct sum decompositions as
$R(T)^W$ modules: \be\label{4} R(T)_{\mathfrak{p}}^{W_{\Delta\setminus
    I}}=\bigoplus_{J\subseteq I} R(T)_{J}\ee \be\label{5}
R(T)_{\mathfrak{p}}^{W_{\Delta\setminus I}}=(\sum_{J\subsetneq I}
R(T)_{\mathfrak{p}}^{W_{\Delta \setminus J}})\bigoplus R(T)_{I} \ee
for $I\subseteq \Delta$.\elem {\bf Proof:} By (\ref{1}) we note that
$W^{\Delta\setminus I}=\bigsqcup_{J\subseteq I} C^J$. Hence Lemma
\ref{basis2} implies that:
\be\label{3}R(T)_{\mathfrak{p}}^{W_{\Delta\setminus
    I}}=\bigoplus_{J\subseteq I}\bigoplus_{v\in
  C^J}R(T)_{\mathfrak{p}}^{W}\cdot r_v.\ee Using (\ref{2}), the
proof of the lemma is now exactly same as that of Lemma 1.10 of
\cite{u}. $\Box$

In $R(T)_{\mathfrak{p}}$ we have: \be \label{multr}r_{v}\cdot
r_{v^{\prime}}=\sum_{J\subseteq (I\cup I^{\prime})} \sum_{w\in
  C^{J}}a^{w}_{v,v^{\prime}}\cdot r_{w}\ee for certain elements
$a^{w}_{v,v^{\prime}}\in
R(G)_{\mathfrak{p}}=R(T)_{\mathfrak{p}}^{W}$,
$~\forall~ v\in C^{I}$, $v^{\prime}\in C^{I^{\prime}}$ and
$w\in C^{J}$ $J\subseteq (I\cup I^{\prime})$.

Further, let $\lambda_I$ denote the image of the element
$\prod_{\alpha\in I}(1-e^{-{\alpha}})\in R(T)$ in $K(G/B)$
under$c_K$
for every $I\subseteq \Delta$. 

Further, let
\[K(G/B)_{\bq,I}:=\bigoplus_{v\in C^{I}} \bq \cdot
[{\co_{X^v}}].\]
Then we have:
\[K(G/B)_{\bq}=\bigoplus_{I\subseteq \Delta}K(G/B)_{I}.\]

\subsection{Main Theorem}
Let $X:=\bar{G_{ad}}$ denote the wonderful compactification of
the
semisimple adjoint group $G_{ad}=G/Z(G)$, where $Z(G)$ denotes
the
center of $G$, constructed by De Concini and Procesi in
\cite{dp}.

Note that $K_{G\times G}(X)$ is an $\cR:=R(G)\otimes
R(G)$-module. Let
\[\cs:=R(G)\otimes R(G)_{\mathfrak{p}}.\] Then we note that the
forgetful
homomorphism extends to \be\label{extchar1}f: K_{G\times
  G}(X)\otimes_{\cR}\cs\ra K(X)_{\bq}.\ee

\bth\label{kdec} The ring $K_{G\times G}(X)\otimes_{\cR}\cs$ has the
following direct sum decomposition as an $\cs$-module: \be K_{G\times
  G}(X)\otimes_{\cR}\cs=\bigoplus_{I\subseteq \Delta}\prod_{\alpha\in
  I} (1-e^{\alpha(u)})\cdot R(T)\otimes R(T)_{I}.\ee Further, the
above direct sum is a free $R(T)\otimes R(G)_{\mathfrak{p}}$-module of
rank $|W|$ with basis \[\{\prod_{\alpha\in I}(1-e^{\alpha(u)})\otimes
r_{v} :~v\in C^{I}~and~ I\subseteq \Delta\},\] where $C^{I}$ is as
defined in (\ref{1}) and $\{r_{v}\}$ is as defined above.  Moreover,
we can identify the component $R(T)\otimes 1\subseteq R(T)\otimes
R(G)_{\mathfrak{p}}$ in the above direct sum with the subring of
$K_{G\times G}(X)$ generated by $Pic^{G\times G}(X)$. (We refer to
\cite{str} for similar description of the equivariant cohomology ring
of the wonderful compactifications.)\eeth

{\bf Proof:} Recall from Lemma 3.2 of \cite{u} that we have a chain of
inclusions: \be\label{wondch} R(T)\otimes R(G)\subseteq K_{G\times
  G}(X)\subseteq R(T)\otimes R(T)\ee where $K_{G\times G}(X)$ consists
of elements $f(u,v)\in R(T)\otimes R(T)$ that satisfy the relations
\be\label{wondrel}(1,s_{\alpha})f(u,v)\equiv
f(u,v)~(\mbox{mod}(1-e^{\alpha(u)})) ~\mbox{for~
  every~}\alpha\in\Delta.\ee Now, since $1\otimes R(G)_{\mathfrak{p}}$
is a flat $1\otimes R(G)$-module, we see that $\cs$ is flat as an
$\cR$-module.  This implies from
(\ref{wondch}) that we have a chain of inclusions:
\be\label{wondchloc} R(T)\otimes R(G)\otimes_{\cR}\cs\subseteq
K_{G\times G}(X)\otimes_{\cR}\cs\subseteq R(T)\otimes
R(T)\otimes_{\cR}\cs.\ee Further, $K_{G\times G}(X)\otimes_{\cR}\cs$
consists of elements \[f(u,v)\in R(T)\otimes R(T)\otimes_{\cR}\cs\]
that satisfy the relations
\be\label{wondrelloc}(1,s_{\alpha})f(u,v)\equiv
f(u,v)~(\mbox{mod}(1-e^{\alpha(u)})) ~\mbox{for~
  every~}\alpha\in\Delta.\ee Here we use the fact that if $f(u,v)\in
\cs$, then $(1,s_{\alpha})f(u,v)=f(u,v)$ for every $\alpha\in\Delta$.
The theorem now follows by using Lemma \ref{dirsum} above, and
replacing the Steinberg basis $\{f_v\}_{v\in W^I}$ by the canonical
lifting of the Schubert basis $\{r_v\}_{v\in W^{I}}$, in the proof of
Theorem 3.8 of \cite{u}.  $\Box$

\bth\label{main} The subring of $K(X)$ generated by classes of
line
bundles is isomorphic to $K(G/B)$. Moreover, $K(X)$ is a free
module
of rank $|W|$ over $K(G/B)$.  More explicitly, let
\[\gamma_{v}:=1\otimes [\co_{X^v}]\in K(G/B)\otimes
K(G/B)_{\bq,I}\] for $v\in C^{I}$ for every $I\subseteq
\Delta$.Then
we have:
\[K(X)_{\bq}\simeq \bigoplus_{v\in W} K(G/B)\cdot \gamma_v .\]
Further, the above isomorphism is a ring isomorphism, where the
multiplication of any two basis elements $\gamma_v$ and
$\gamma_{v'}$
is defined as follows:
\[\gamma_v\cdot \gamma_{v'}:=\sum_{J\subseteq (I\cup
I^{\prime})}\sum_{w\in C^{J}}(\lambda_{I\cap I'}\cdot
\lambda_{{(I\cup
I')\setminus J}}\cdot c^{w}_{v,v'}) \cdot \gamma_{w}.\]
where $c^w_{v,v^{\prime}}\in\bz$ are as defined in
(\ref{ostrconst1}).\eeth
        
{\bf Proof:} Since $c_K(r_v)=[\co_{X^{v}}]$ for $v\in W^{I}$
and$I\subseteq \Delta$, we note that the image under $c_K$ of
the
element
$a^{w}_{v,v^{\prime}}\in R(G)_{\mathfrak{p}}$ defined in
(\ref{multr})
is nothing but the structure constant
$c^w_{v,v^{\prime}}\in\bz$ defined in (\ref{ostrconst1}). The
proof now follows exactly as
that
of Theorem 3.12 on p.403 of \cite{u}. $\Box$

\brem\label{comp} Note that Theorem \ref{main} is a restatement
of
Theorem 3.12 of \cite{u}, obtained by replacing the Steinberg
basis
$\{f_v\}_{v\in W^I}$ by the lift of the Schubert basis
$\{r_v\}_{v\in
W^{I}}$. In Theorem 3.12 of \cite{u}, the multiplicative
structure
of the ordinary $K$-ring of the wonderful group
compactifications was
described in terms of the structure constants of the image of
the
Steinberg basis $\{f_v\}_{v\in W^I}$ under $c_K$. These
structure
constants do not have any known relations to geometry or
representation theory.

Whereas now we see that the multiplicative structure constants
of the
basis $\gamma_v=1\otimes [\co_{X^v}]$ of $K(X)_{\bq}$ as
$K(G/B)$-module are determined explicitly in terms of the
multiplicative structure constants of the Schubert basis
$c_{K}(r_v)$
described above in Prop. {\ref{oMultschubert}}. These structure
constants have been described in \S2 and \S3 above, and also
are known
to have nice geometric and representation theoretic
interpretations
(see for example \cite{br3} and \cite{ls}).\erem

\section{Appendix}

\subsection{An Explicit lift of $[\co_{X^{w_0}}]_{T}$ in $R(T)\otimes
  R(T)$}

Let $\{e^{p_w}\}_{w\in W}$ be the basis defined by Steinberg of
$R(T)$
as an $R(T)^{W}$-module where:
$$p_w=w(\sum_{{\alpha\in \Delta},
w(\alpha)<0}\omega_{\alpha})$$for
$w\in W$ (see \cite{st}). Then by lemme 4 and prop. 3 of
\cite{m} we
see that the matrix $M=(L_{w'}(e^{p_w}))_{w,w'\in W}$ with
entries in
$R(T)$ is invertible. Thus there exists a unique vector
$(a_w)_{w\in
  W}$ such that
\be\label{explift}\sum_{w\in W} a_{w}\cdot
L_{w'}(e^{p_w})=e^{-\rho}\cdot
\delta_{w',w_0}\ee for every $w'\in W$.  Now, defining
$b_w:=e^{\rho+p_w}$, we see that the element
$$v_0=\sum_{w\in W}a_w\otimes b_w$$ in $R(T)\otimes R(T)$
satisfies (\ref{eq3}). Thus we have a canonical choice of an element
\be\label{explift1} u_0=v_0\cdot (1\otimes e^{-\rho})=\sum_{w\in
  W}a_w\otimes e^{p_w}\ee in $R(T)\otimes R(T)$ which satisfies
Prop. \ref{slift}. We now illustrate the computation of
$u_0=\sum_{w\in W}a_w\otimes b_w$ in $R(T)\otimes R(T)$ for the case
when $G$ is of type $A_2$. \beg\label{example1lift} When $G$ is of
type $A_2$, we have $\Delta=\{\alpha,\beta\}$ and $\omega_{\alpha}$
and $\omega_{\beta}$ are the fundamental weights dual to $\alpha$ and
$\beta$ respectively. Further, $W=\{1,s_\alpha,
s_{\beta},s_{\alpha}s_{\beta},
s_{\beta}s_{\alpha},s_{\alpha}s_{\beta}s_{\alpha}\}$, where
$s_{\alpha}$ and $s_{\beta}$ are the simple reflections corresponding
to $\alpha$ and $\beta$ satisfying the braid relation
$(s_{\alpha}s_{\beta})^3=1$. Moreover, we note that
$\rho=\omega_{\alpha}+\omega_{\beta}$. In this case the Steinberg
basis are determined to be:
\be\label{st}\begin{array}{ccccccc}&e^{p_1}&=&1\\
  &e^{p_{s_{\alpha}}}&=&e^{\omega_{\beta}-\omega_{\alpha}}\\
  &e^{p_{s_{\beta}}}&
  =&e^{\omega_{\alpha}-\omega_{\beta}}\\
  &e^{p_{s_{\alpha}s_{\beta}}}&=&e^{-{\omega_{\alpha}}}
  \\ &e^{p_{s_{\beta}s_{\alpha}}}&=&e^{-\omega_{\beta}}\\
  &e^{p_{s_{\alpha}s_{\beta}s_{\alpha}}}&=&e^{-\omega_{\alpha}-\omega_{\beta}}.
 \end{array}\ee Furthermore,  the matrix   
 $(L_{w'}(e^{p_{w}}))$ is determined to be: 
 \be\label{matrix}\left(\begin{array}{cccccccccc}1&
e^{\omega_{\beta}-\omega_{\alpha}}&
e^{\omega_{\alpha}-\omega_{\beta}}&
e^{-\omega_{\alpha}}& e^{-\omega_{\beta}}&
e^{-\omega_{\alpha}-\omega_{\beta}}\\ 0 &
e^{\omega_{\beta}-\omega_{\alpha}}&
-e^{-\omega_{\alpha}}&e^{-\omega_{\alpha}}&0&
e^{-\omega_{\alpha}-\omega_{\beta}}\\0 &-e^{-\omega_{\beta}}&
e^{\omega_{\alpha}-\omega_{\beta}}&0 &e^{-\omega_{\beta}}&
e^{-\omega_{\alpha}-\omega_{\beta}}\\ 0 & 0 &-
e^{-\omega_{\alpha}}&0 &0&
e^{-\omega_{\alpha}-\omega_{\beta}}\\ 0 & -e^{-\omega_{\beta}}&
0 &0 &0&
     e^{-\omega_{\alpha}-\omega_{\beta}}\\ 0 & 0 &0& 0&0&
e^{-\omega_{\alpha}-\omega_{\beta}}\end{array}\right)\ee where
the
rows of the matrix correspond respectively to $L_{w'}(e^{p_w})$
for
 $w'\in W$ for the above ordering. 
We can now solve the system of equations (\ref{explift}) to get
\be\label{a_w's}\begin{array}{cccc} &a_1
&=&-e^{-\omega_{\alpha}-\omega_{\beta}} \\
&a_{s_{\alpha}}&=&e^{-\omega_{\alpha}} \\
&a_{s_{\beta}}&=&e^{-\omega_{\beta}}
\\&a_{s_{\alpha}s_{\beta}}&=&-e^{-\omega_{\alpha}+\omega_{\beta}}\\
&a_{s_{\beta}s_{\alpha}}&=&-e^{-\omega_{\beta}+\omega_{\alpha}}\\
   &a_{s_{\alpha}s_{\beta}s_{\alpha}}&=&1
 \end{array}\ee Hence from (\ref{explift1}) we get 
 \be\label{expliftA_2}\begin{array}{lllll}u_0=&1\otimes 
   e^{-\omega_{\alpha}-\omega_{\beta}}&- 
   &e^{-\omega_{\alpha}-\omega_{\beta}}\otimes 1 &+\\ &  
e^{-\omega_{\alpha}}\otimes
e^{\omega_{\beta}-\omega_{\alpha}}&-&e^{\omega_{\beta}-\omega_{\alpha}}\otimes
e^{-\omega_{\alpha}}&+\\&
e^{-\omega_{\beta}}\otimes
e^{\omega_{\alpha}-\omega_{\beta}}&-&e^{\omega_{\alpha}-\omega_{\beta}}\otimes
e^{-\omega_{\beta}}. \end{array}\ee\eeg

\subsection{Examples for computations of structure constants}
We now illustrate the computations in Lemma \ref{Multschubert}
and in
Lemma \ref{chevalley} respectively by the following examples,
when
 $G$ is of the type $A_2$. We follow the notations of Example
\ref{example1lift}. \beg\label{example1} Let $u_0=\sum_{w\in W}
a_w\otimes e^{p_w}\in R(T)\otimes R(T)$ be the lift of
 $[\co_{X^{w_0}}]_{T}$ as in (\ref{expliftA_2}). Further, let
 $x=s_{\alpha}$, $y=s_{\alpha}s_{\beta}$, and
\[t^{w,w'}_{x,y}:=L_{s_{\beta}s_{\alpha}}(b_w\cdot
e^{-\rho})\cdot
 L_{s_{\alpha}}(b_{w'}\cdot e^{-\rho})\cdot
 e^{\rho}=L_{s_{\beta}s_{\alpha}}(e^{p_w})\cdot
 L_{s_{\alpha}}(e^{p_{w'}}) \cdot e^{\rho}.\] Then by
(\ref{strconst1eq}), the multiplicative structures constants of
$[\co_{x}]_{T}\cdot[\co_{y}]_{T}$ in $K_{T}(G/B)$ are obtained
 recursively as follows: \be\label{str}\begin{array}{lllllll}
   C^1_{x,y}&= \sum_{w,w'}a_w\cdot
   a_{w'} \cdot e^{\rho}\cdot t^{w,w'}_{x,y},\\
   C^{s_{\alpha}}_{x,y}&=\sum_{w,w'}a_w\cdot a_{w'} \cdot
   e^{\rho}\cdot L_{s_{\alpha}}(t^{w,w'}_{x,y})-C^1_{x,y},\\
C^{s_{\beta}}_{x,y}&=\sum_{w,w'}a_w\cdot a_{w'} \cdot
e^{\rho}\cdot
   L_{s_{\beta}}(t^{w,w'}_{x,y})-C^1_{x,y},\\
C_{x,y}^{s_{\alpha}s_{\beta}}&=\sum_{w,w'}a_w\cdot a_{w'} \cdot
e^{\rho}\cdot
L_{s_{\alpha}s_{\beta}}(t^{w,w'}_{x,y})-C_{x,y}^{s_{\alpha}}-C_{x,y}^{s_{\beta}}-C^1_{x,y},\\
C_{x,y}^{s_{\beta}s_{\alpha}}&=\sum_{w,w'}a_w\cdot a_{w'} \cdot
e^{\rho}\cdot
L_{s_{\beta}s_{\alpha}}(t^{w,w'}_{x,y})-C_{x,y}^{s_{\alpha}}-C_{x,y}^{s_{\beta}}-C^1_{x,y},\\
C_{x,y}^{s_{\alpha}s_{\beta}s_{\alpha}}&=\sum_{w,w'}a_w\cdot
a_{w'}
   \cdot e^{\rho}\cdot
L_{s_{\alpha}s_{\beta}s_{\alpha}}(t^{w,w'}_{x,y})-C^{s_{\alpha}s_{\beta}}_{x,y}-C^{s_{\beta}s_{\alpha}}_{x,y}+
\\ &~~~~
-C_{x,y}^{s_{\alpha}}-C_{x,y}^{s_{\beta}}-C^1_{x,y}.\end{array}\ee
Now, from (\ref{matrix}) it follows that
$t_{x,y}^{w,w'}=0$ when
$w=1, s_{\beta}, s_{\alpha}s_{\beta}, s_{\beta}s_{\alpha}$ or
$w'=1,s_{\beta}s_{\alpha}$.  Further, we have:
\be\label{t}\begin{array}{lllll}t_{x,y}^{s_{\alpha},s_{\alpha}}&=&-e^{\omega_{\beta}}\\t_{x,y}^{s_{\alpha},s_{\beta}}&=&
  1\\
  t_{x,y}^{s_{\alpha},s_{\alpha}s_{\beta}}&=&
  -1\\t_{x,y}^{s_{\alpha},
    s_{\alpha}s_{\beta}s_{\alpha}}&=&-e^{-\omega_{\beta}}\\
  t_{x,y}^{s_{\alpha}s_{\beta}s_{\alpha},s_{\alpha}}&
  =&e^{\omega_{\beta}-\omega_{\alpha}}\\
t_{x,y}^{s_{\alpha}s_{\beta}s_{\alpha},s_{\beta}}&=&-e^{-\omega_{\alpha}}\\
  t_{x,y}^{s_{\alpha}s_{\beta}s_{\alpha},s_{\alpha}s_{\beta}}&
=&e^{-\omega_{\alpha}}\\t_{x,y}^{s_{\alpha}s_{\beta}s_{\alpha},s_{\alpha}s_{\beta}s_{\alpha}}&=&e^{-\omega_{\alpha}-\omega_{\beta}}.
\end{array}\ee

\be\label{s{alpha}t}\begin{array}{lllll}L_{s_{\alpha}}(t_{x,y}^{s_{\alpha},s_{\alpha}})&
  =&0\\L_{s_{\alpha}}(t_{x,y}^{s_{\alpha},s_{\beta}})&=&
  0\\
  L_{s_{\alpha}}(t_{x,y}^{s_{\alpha},s_{\alpha}s_{\beta}})&=&
  0\\L_{s_{\alpha}}(t_{x,y}^{s_{\alpha},
    s_{\alpha}s_{\beta}s_{\alpha}})&=&0\\
L_{s_{\alpha}}(t_{x,y}^{s_{\alpha}s_{\beta}s_{\alpha},s_{\alpha}})&
  =&e^{\omega_{\beta}-\omega_{\alpha}}\\
L_{s_{\alpha}}(t_{x,y}^{s_{\alpha}s_{\beta}s_{\alpha},s_{\beta}})&=
  &-e^{-\omega_{\alpha}}\\
L_{s_{\alpha}}(t_{x,y}^{s_{\alpha}s_{\beta}s_{\alpha},s_{\alpha}s_{\beta}})&
=&e^{-\omega_{\alpha}}\\L_{s_{\alpha}}(t_{x,y}^{s_{\alpha}s_{\beta}s_{\alpha},
s_{\alpha}s_{\beta}s_{\alpha}})&=&e^{-\omega_{\alpha}-\omega_{\beta}}.
\end{array}\ee

\be\label{s{beta}t}\begin{array}{lllll}
  L_{s_{\beta}}(t_{x,y}^{s_{\alpha},s_{\alpha}})&
  =&e^{-\omega_{\beta}+\omega_{\alpha}}\\
  L_{s_{\beta}}(t_{x,y}^{s_{\alpha},s_{\beta}})&=&
  0\\
  L_{s_{\beta}}(t_{x,y}^{s_{\alpha},s_{\alpha}s_{\beta}})&=&
  0\\L_{s_{\beta}}(t_{x,y}^{s_{\alpha},
    s_{\alpha}s_{\beta}s_{\alpha}})&=&-e^{-\omega_{\beta}}\\
L_{s_{\beta}}(t_{x,y}^{s_{\alpha}s_{\beta}s_{\alpha},s_{\alpha}})&
  =&-e^{-\omega_{\beta}}\\
L_{s_{\beta}}(t_{x,y}^{s_{\alpha}s_{\beta}s_{\alpha},s_{\beta}})&=&0\\
L_{s_{\beta}}(t_{x,y}^{s_{\alpha}s_{\beta}s_{\alpha},s_{\alpha}s_{\beta}})&
  =&0\\L_{s_{\beta}}(t_{x,y}^{s_{\alpha}s_{\beta}s_{\alpha},
s_{\alpha}s_{\beta}s_{\alpha}})&=&e^{-\omega_{\alpha}-\omega_{\beta}}.
\end{array}\ee

\be\label{s{alpha}s{beta}t}\begin{array}{lllllll}
  L_{s_{\alpha}s_{\beta}}(t_{x,y}^{s_{\alpha},s_{\alpha}})&
  =& -e^{-\omega_{\alpha}}\\
  L_{s_{\alpha}s_{\beta}}(t_{x,y}^{s_{\alpha},s_{\beta}})&=&0\\
L_{s_{\alpha}s_{\beta}}(t_{x,y}^{s_{\alpha},s_{\alpha}s_{\beta}})&=&
  0\\L_{s_{\alpha}s_{\beta}}(t_{x,y}^{s_{\alpha},
      s_{\alpha}s_{\beta}s_{\alpha}})&=&0\\
L_{s_{\alpha}s_{\beta}}(t_{x,y}^{s_{\alpha}s_{\beta}s_{\alpha},s_{\alpha}})&
    =&0\\
L_{s_{\alpha}s_{\beta}}(t_{x,y}^{s_{\alpha}s_{\beta}s_{\alpha},s_{\beta}})&=&0\\
L_{s_{\alpha}s_{\beta}}(t_{x,y}^{s_{\alpha}s_{\beta}s_{\alpha},s_{\alpha}s_{\beta}})&
=&0\\L_{s_{\alpha}s_{\beta}}(t_{x,y}^{s_{\alpha}s_{\beta}s_{\alpha},
s_{\alpha}s_{\beta}s_{\alpha}})&=&e^{-\omega_{\alpha}-\omega_{\beta}}.
\end{array}\ee

\be\label{s{beta}s{alpha}t}\begin{array}{llll}
  L_{s_{\beta}s_{\alpha}}(t_{x,y}^{s_{\alpha},s_{\alpha}})&
  =0\\
  L_{s_{\beta}s_{\alpha}}(t_{x,y}^{s_{\alpha},s_{\beta}})&=
 0\\
L_{s_{\beta}s_{\alpha}}(t_{x,y}^{s_{\alpha},s_{\alpha}s_{\beta}})&=0\\
  L_{s_{\beta}s_{\alpha}}(t_{x,y}^{s_{\alpha},
    s_{\alpha}s_{\beta}s_{\alpha}})&=0\\
L_{s_{\beta}s_{\alpha}}(t_{x,y}^{s_{\alpha}s_{\beta}s_{\alpha},s_{\alpha}})&
  =-e^{-\omega_{\beta}}\\
L_{s_{\beta}s_{\alpha}}(t_{x,y}^{s_{\alpha}s_{\beta}s_{\alpha},s_{\beta}})&=0\\
L_{s_{\beta}s_{\alpha}}(t_{x,y}^{s_{\alpha}s_{\beta}s_{\alpha},s_{\alpha}s_{\beta}})&
  =0\\
L_{s_{\beta}s_{\alpha}}(t_{x,y}^{s_{\alpha}s_{\beta}s_{\alpha},s_{\alpha}s_{\beta}s_{\alpha}})&=e^{-\omega_{\alpha}-\omega_{\beta}}
\end{array}\ee

\be\label{s{alpha}s{beta}s{alpha}t}\begin{array}{llll}
&L_{s_{\alpha}s_{\beta}s_{\alpha}}(t_{x,y}^{s_{\alpha},s_{\alpha}})&
  =0\\
&L_{s_{\alpha}s_{\beta}s_{\alpha}}(t_{x,y}^{s_{\alpha},s_{\beta}})&=
  0\\
&L_{s_{\alpha}s_{\beta}s_{\alpha}}(t_{x,y}^{s_{\alpha},s_{\alpha}s_{\beta}})&=0\\
  &L_{s_{\alpha}s_{\beta}s_{\alpha}}(t_{x,y}^{s_{\alpha},
    s_{\alpha}s_{\beta}s_{\alpha}})&=0\\
&L_{s_{\alpha}s_{\beta}s_{\alpha}}(t_{x,y}^{s_{\alpha}s_{\beta}s_{\alpha},s_{\alpha}})&
  =0\\
&L_{s_{\alpha}s_{\beta}s_{\alpha}}(t_{x,y}^{s_{\alpha}s_{\beta}s_{\alpha},s_{\beta}})&=0\\
&L_{s_{\alpha}s_{\beta}s_{\alpha}}(t_{x,y}^{s_{\alpha}s_{\beta}s_{\alpha},s_{\alpha}s_{\beta}})&
  =0\\
&L_{s_{\alpha}s_{\beta}s_{\alpha}}(t_{x,y}^{s_{\alpha}s_{\beta}s_{\alpha},
s_{\alpha}s_{\beta}s_{\alpha}})&=e^{-\omega_{\alpha}-\omega_{\beta}}
\end{array}\ee

\be\label{multa_w's}\begin{array}{lllll}e^{\rho}\cdot
  a_{s_{\alpha}}\cdot a_{s_{\alpha}}&
  =&e^{\omega_{\beta}-\omega_{\alpha}}\\e^{\rho}\cdot
  a_{s_{\alpha}}\cdot a_{s_{\beta}}&=&
  1\\
  e^{\rho}\cdot a_{s_{\alpha}}\cdot a_{s_{\alpha}s_{\beta}}&=&
-e^{2\omega_{\beta}-\omega_{\alpha}}\\e^{\rho}\cdot
a_{s_{\alpha}}\cdot
  a_{s_{\alpha}s_{\beta}s_{\alpha}}&=&e^{\omega_{\beta}}\\
e^{\rho}\cdot a_{s_{\alpha}s_{\beta}s_{\alpha}}\cdot
a_{s_{\alpha}}&
=&e^{\omega_{\beta}}\\e^{\rho}\cdot
a_{s_{\alpha}s_{\beta}s_{\alpha}}\cdot
a_{s_{\beta}}&=&e^{\omega_{\alpha}}\\
e^{\rho}\cdot a_{s_{\alpha}s_{\beta}s_{\alpha}}\cdot
a_{s_{\alpha}s_{\beta}}& =&-e^{2\omega_{\beta}}\\e^{\rho}\cdot
a_{s_{\alpha}s_{\beta}s_{\alpha}}\cdot
a_{s_{\alpha}s_{\beta}s_{\alpha}}&=&e^{\rho}.
\end{array}\ee
Now, substituting (\ref{t}) to (\ref{multa_w's}) in (\ref{str})
we get: \be\label{str1} \begin{array}{llll} C^1_{x,y}&= 0\\
  C^{s_{\alpha}}_{x,y}&=0\\
  C^{s_{\beta}}_{x,y}&=0\\
  C_{x,y}^{s_{\alpha}s_{\beta}}&=1-e^{-\alpha}\\
  C_{x,y}^{s_{\beta}s_{\alpha}}&=0\\
C_{x,y}^{s_{\alpha}s_{\beta}s_{\alpha}}&=1-(1-e^{-\alpha})=e^{-\alpha}.
\end{array}\ee
\eeg 

\brem\label{compgriffethram} In particular, note that:
\[(-1)^{l(x)+l(y)+l(s_{\alpha}s_{\beta})}\cdot
C_{x,y}^{s_{\alpha}s_{\beta}}=e^{-\alpha}-1\] and
\[(-1)^{l(x)+l(y)+l(s_{\alpha}s_{\beta}s_{\alpha})}\cdot
C_{x,y}^{s_{\alpha}s_{\beta}s_{\alpha}}=(e^{-\alpha}-1)+1.\]
This
verifies Conjecture 3.10 in \cite{gk} (conjecture of
Griffeth-Ram) for
this example. We refer the reader to \S5 of \cite{gr} for the
computation of multiplicative structure constants in rank $2$
cases
using different methods, and a proof of the positivity
conjecture in
\cite{agm}.\erem

\beg\label{example2} Let $x=w_0=s_{\alpha}s_{\beta}s_{\alpha}$
and
$\lambda=\rho$. Let $u_0$ be the lift of $[\co_{X^{w_0}}]_{T}$
as in
(\ref{expliftA_2}). In particular, we note that
$e^{\lambda}\cdot
L_{x^{-1}w_0}(b_w\cdot e^{-\rho})=b_w$. Then the Chevalley
structure
constants of $[\co_{X^x}]_{T}\cdot [\cl^{T}(\rho)]_{T}$ are
obtained
recursively as
follows: 

\be\label{chevst}\begin{array}{llll}Q^{\rho}_{x,1}&=\sum_{w\in
    W}e^{\rho}\cdot
  a_w\cdot b_w\\
  Q^{\rho}_{x,s_{\alpha}}&=\sum_{w\in W} e^{\rho}\cdot a_w\cdot
L_{s_{\alpha}}(b_w)-Q^{\rho}_{x,1}\\Q^{\rho}_{x,s_{\beta}}&=\sum_{w\in
    W}e^{\rho}\cdot a_w\cdot
L_{s_{\beta}}(b_w)-Q^{\rho}_{x,1}\\Q^{\rho}_{x,s_{\alpha}s_{\beta}}&=\sum_{w\in
    W}e^{\rho}\cdot
a_w\cdot
L_{s_{\alpha}s_{\beta}}(b_w)-Q^{\rho}_{x,s_{\alpha}}-Q^{\rho}_{x,s_{\beta}}-Q^{\rho}_{x,1}\\
  Q^{\rho}_{x,s_{\beta}s_{\alpha}}&=\sum_{w\in W}e^{\rho}\cdot
  a_w\cdot
L_{s_{\beta}s_{\alpha}}(b_w)-Q^{\rho}_{x,s_{\alpha}}-Q^{\rho}_{x,s_{\beta}}-Q^{\rho}_{x,1}\\
  Q^{\rho}_{x,s_{\alpha}s_{\beta}s_{\alpha}}&=\sum_{w\in
    W}e^{\rho}\cdot a_w\cdot
L_{s_{\alpha}s_{\beta}s_{\alpha}}(b_w)-Q^{\rho}_{x,s_{\alpha}s_{\beta}}-Q^{\rho}_{x,s_{\beta}s_{\alpha}}-
  Q^{\rho}_{x,s_{\alpha}}+\\&~~~~
-Q^{\rho}_{x,s_{\beta}}-Q^{\rho}_{x,1}.\end{array}\ee Now,
since$b_w=e^{\rho+p_w}$, from (\ref{st}) it follows that:
\be\label{b_w's}\begin{array}{cccc} &b_1
  &=&e^{\rho}\\ &b_{s_{\alpha}}&=&
e^{2\omega_{\beta}} \\ &b_{s_{\beta}}&=&e^{2\omega_{\alpha}} \\
&b_{s_{\alpha}s_{\beta}}&=&e^{\omega_{\beta}}\\
  &b_{s_{\beta}s_{\alpha}}&=&e^{\omega_{\alpha}}\\
  &b_{s_{\alpha}s_{\beta}s_{\alpha}}&=&1
\end{array}\ee Hence by substituting (\ref{a_w's}) and
(\ref{b_w's}) in
(\ref{chevst}) we get:
\be\label{chevst1}\begin{array}{llll}Q^{\rho}_{x,1}&=0\\
Q^{\rho}_{x,s_{\alpha}}&=0\\Q^{\rho}_{x,s_{\beta}}&=0\\Q^{\rho}_{x,s_{\alpha}s_{\beta}}&=0\\
  Q^{\rho}_{x,s_{\beta}s_{\alpha}}&=0\\
Q^{\rho}_{x,s_{\alpha}s_{\beta}s_{\alpha}}&=e^{-2\rho}.\end{array}\ee
\eeg

\brem\label{compkk}It was recently brought to the notice of the
author
that some results in this article coincide with or follow from
results
in the paper by Kostant and Kumar \cite{kk}. We wish to mention
that
these results have been independently proved during the course
of this
work using slightly different techniques. We give here the
exact cross
references for the benefit of the reader. Here Lemma 2.1 is
essentially Theorem 4.4 of \cite{kk}. Also, Prop. 2.5 follows
essentially from Lemma 4.12 of \cite{kk} and Lemma 2.9 is the
analogue
of Prop. 2.25 of \cite{kk} where the structure constants are
determined with respect to the dual of the structure sheaf
basis.

\erem

\end{document}